\documentclass[12pt]{article}
\newcommand{\ba}{\begin{array}}
\newcommand{\ea}{\end{array}}
\newcommand{\be}{\begin{equation}}
\newcommand{\ee}{\end{equation}}
\newcommand{\la}{\label}
\newcommand{\bea}{\begin{eqnarray}}
\newcommand{\eea}{\end{eqnarray}}
\newcommand{\ch}{\choose}
\newcommand{\n}{\nonumber}
\newcommand{\nn}{\nonumber \\}
\newcommand{\ds}{\displaystyle}
\newcommand{\ndots}{n=0,1,2,\ldots}
\renewcommand{\a}{\alpha}
\renewcommand{\b}{\beta}
\newcommand{\G}{\Gamma}

\newcommand{\GP}{P_n^{\a,\b,M,N}(x)}
\newcommand{\SP}{P_n^{(\a,\a)}(x)}
\newcommand{\SGP}{P_n^{\a,\a,M,M}(x)}
\renewcommand{\l}{\left}
\renewcommand{\r}{\right}
\newcommand{\set}[1]{\left\{#1\right\}_{n=0}^{\infty}}
\newcommand{\hyp}[5]{\mbox{}_{#1}F_{#2}
\left(\left.{#3 \atop #4}\right|#5\right)}

\begin{document}

\begin{center}
{ \Large \bf Finding differential equations for symmetric generalized
ultraspherical polynomials by using inversion methods }
\end{center}
\begin{center}
{ \Large \bf J.~Koekoek and R.~Koekoek }
\end{center}

\vspace*{1cm}

\begin{abstract}
We find all differential equations of the form
$$M\sum_{i=0}^{\infty}a_i(x)y^{(i)}(x)+
(1-x^2)y''(x)-2(\a+1)xy'(x)+n(n+2\a+1)y(x)=0,$$
where the coefficients $\l\{a_i(x)\r\}_{i=1}^{\infty}$ are independent of
$n$ and $a_0(x):=a_0(n,\a)$ is independent of $x$, satisfied by the
symmetric generalized ultraspherical polynomials $\set{\SGP}$ which are
orthogonal on the interval $[-1,1]$ with respect to the weight function
$$\frac{\G(2\a+2)}{2^{2\a+1}\l\{\G(\a+1)\r\}^2}\l(1-x^2\r)^{\a}+
M\l[\delta(x+1)+\delta(x-1)\r],$$
where $\a>-1$ and $M\ge 0$.

In order to find explicit formulas for the coefficients of these
differential equations we have to solve systems of equations of the form
$$\sum_{i=1}^{\infty}A_i(x)D^i\SP=F_n(x),\;n=1,2,3,\ldots,$$
where the coefficients $\l\{A_i(x)\r\}_{i=1}^{\infty}$ are independent of
$n$. This system of equations has a unique solution given by
$$A_i(x)=2^i\sum_{j=1}^i\frac{2\a+2j+1}{(2\a+j+1)_{i+1}}
P_{i-j}^{(-\a-i-1,-\a-i-1)}(x)F_j(x),\;i=1,2,3,\ldots.$$
This is a consequence of the inversion formula
\bea & &\sum_{k=j}^i\frac{2\a+2k+1}{(2\a+k+j+1)_{i-j+1}}\times{}\nn
& &{}\hspace{1cm}{}\times
P_{i-k}^{(-\a-i-1,-\a-i-1)}(x)P_{k-j}^{(\a+j,\a+j)}(x)=\delta_{ij},
\;j\le i,\;i,j=0,1,2,\ldots.\n\eea
\end{abstract}

\vspace*{1cm}


\section{Introduction}

Let $\a>-1$. In \cite{Symjac} we found all differential equations of the
form
\be\la{DVSGP}M\sum_{i=0}^{\infty}a_i(x)y^{(i)}(x)+
(1-x^2)y''(x)-2(\a+1)xy'(x)+n(n+2\a+1)y(x)=0,\ee
where the coefficients $\l\{a_i(x)\r\}_{i=0}^{\infty}$ are continuous
functions on the real line and $\l\{a_i(x)\r\}_{i=1}^{\infty}$ are
independent of $n$, satisfied by the symmetric generalized ultraspherical
polynomials $\set{\SGP}$ defined by
$$\SGP=C_0\SP-C_1xD\SP,\;\ndots,$$
where $\ds D=\frac{d}{dx}$ denotes the differentiation operator and
$$\l\{\ba{l}\ds C_0=1+\frac{2Mn}{\a+1}{n+2\a+1 \ch n}
+4M^2{n+2\a+1 \ch n-1}^2\\
\\
\ds C_1=\frac{2M}{2\a+1}{n+2\a \ch n}
+\frac{2M^2}{\a+1}{n+2\a \ch n-1}{n+2\a+1 \ch n}.\ea\r.$$
The case $2\a+1=0$ must be understood by continuity. These polynomials form
a special case ($\b=\a$ and $N=M$) of the generalized Jacobi polynomials
$\set{\GP}$ introduced by T.H.~Koornwinder in \cite{Koorn}.

In \cite{Report} we gave a proof of the Jacobi inversion formula. The special
case $\b=\a$ of this inversion formula reads
\bea\la{inv}& &\sum_{k=j}^i
\frac{2\a+2k+1}{(2\a+k+j+1)_{i-j+1}}\times{}\nn
& &{}\hspace{1cm}{}\times
P_{i-k}^{(-\a-i-1,-\a-i-1)}(x)P_{k-j}^{(\a+j,\a+j)}(x)=\delta_{ij},
\;j\le i,\;i,j=0,1,2,\ldots.\eea
Again, the case $2\a+1=0$ must be understood by continuity. If we apply this
inversion formula to the system of equations
\be\la{sysalg}\sum_{i=1}^{\infty}A_i(x)D^i\SP=F_n(x),\;n=1,2,3,\ldots,\ee
where the coefficients $\l\{A_i(x)\r\}_{i=1}^{\infty}$ are independent of
$n$, then we find
\be\la{solalg}A_i(x)=2^i\sum_{j=1}^i\frac{2\a+2j+1}{(2\a+j+1)_{i+1}}
P_{i-j}^{(-\a-i-1,-\a-i-1)}(x)F_j(x),\;i=1,2,3,\ldots.\ee

This will be used to find all differential equations of the form
(\ref{DVSGP}), where the coefficients $\l\{a_i(x)\r\}_{i=1}^{\infty}$ are
independent of $n$ and $a_0(x):=a_0(n,\a)$ is independent of $x$.

We will also need the formula
\bea\la{spec}& &\sum_{k=j}^i
\frac{2\a+2k+1}{(2\a+k+j+1)_{i-j+1}}\times{}\nn
& &{}\hspace{1cm}{}\times
P_{i-k}^{(-\a-i-1,-\a-i-1)}(-x)P_{k-j}^{(\a+j,\a+j)}(x)=
\frac{x^{i-j}}{(i-j)!},\;j\le i,\;i,j=0,1,2,\ldots,\hspace{5mm}\eea
which is also proved in \cite{Report}. The case $2\a+1=0$ must be understood
by continuity again.

In this paper we will give the main results. For more details the reader is
referred to the report \cite{Report} where complete proofs are given.

\section{The classical ultraspherical polynomials}

In this section we list the definitions and some properties of the classical
ultraspherical polynomials which we will use in this paper. For details the
reader is referred to \cite{Chihara}, \cite{AS}, \cite{Szego} and the report
\cite{Report}.

The classical ultraspherical polynomials $\set{\SP}$ can be defined by
\bea\la{def1}\SP&=&\sum_{k=0}^n\frac{(n+2\a+1)_k}{k!}
\frac{(\a+k+1)_{n-k}}{(n-k)!}\l(\frac{x-1}{2}\r)^k,\;\ndots\\
&=&\la{def2}(-1)^n\sum_{k=0}^n\frac{(-n-k-2\a)_k}{k!}
\frac{(-n-\a)_{n-k}}{(n-k)!}\l(\frac{x-1}{2}\r)^k,\;\ndots\\
&=&\la{def3}2^{-n}\sum_{k=0}^n{n+\a \ch n-k}{n+\a \ch k}
(x-1)^k(x+1)^{n-k},\;\ndots\eea
for all $\a$. For all $n\in\{0,1,2,\ldots\}$ we have
\be\la{diff}D^i\SP=\frac{(n+2\a+1)_i}{2^i}P_{n-i}^{(\a+i,\a+i)}(x),
\;i=0,1,2,\ldots,n.\ee
The ultraspherical polynomials satisfy the symmetry formula
\be\la{sym}P_n^{(\a,\a)}(-x)=(-1)^n\SP,\;\ndots\ee
and the linear second order differential equation
\be\la{dv}(1-x^2)y''(x)-2(\a+1)xy'(x)+n(n+2\a+1)y(x)=0.\ee

Further we list some formulas involving ultraspherical polynomials which we
will need in this paper. For details the reader is referred to the report
\cite{Report}. First of all we have
\be\la{rel1}2xD\SP=2n\SP+(n+\a)P_{n-2}^{(\a+1,\a+1)}(x),
\;n=2,3,4,\ldots.\ee
Further we have (see for instance \cite{Szego})
\bea\la{rel2} & &(n+2\a+1)(n+2\a+2)P_n^{(\a+1,\a+1)}(x)
-(n+\a)(n+\a+1)P_{n-2}^{(\a+1,\a+1)}(x)\nn
& &{}\hspace{1cm}{}=2(n+\a+1)(2n+2\a+1)\SP,\;n=2,3,4,\ldots.\eea
Finally we will need the formula
\be\la{rel3}(\a+1)P_n^{(\a+1,\a+1)}(x)-(n+\a+1)\SP=
\frac{1}{4}(n+\a+1)(1-x^2)P_{n-2}^{(\a+2,\a+2)}(x),\ee
which also holds for $n=2,3,4,\ldots$.

\section{The computation of the coefficients}

Let $\a>-1$. In \cite{Symjac} we found the coefficients
$\l\{a_i(x)\r\}_{i=0}^{\infty}$ of the differential equation (\ref{DVSGP})
for the symmetric generalized ultraspherical polynomials $\set{\SGP}$. In
order to do this we had to solve the following two systems of equations for
the coefficients $\l\{a_i(x)\r\}_{i=0}^{\infty}$~:
\be\la{sys1}\sum_{i=0}^{\infty}a_i(x)D^i\SP=
\frac{4}{2\a+1}{n+2\a \ch n}D^2\SP\ee
and
\be\la{sys2}\sum_{i=0}^{\infty}ia_i(x)D^i\SP+
x\sum_{i=0}^{\infty}a_i(x)D^{i+1}\SP=4{n+2\a+1 \ch n-1}D^2\SP\ee
for $\ndots$, where the coefficients $\l\{a_i(x)\r\}_{i=0}^{\infty}$ are
continuous functions on the real line and $\l\{a_i(x)\r\}_{i=1}^{\infty}$
are independent of $n$. Now we suppose that $a_0(x):=a_0(n,\a)$ is
independent of $x$ as we did in \cite{Soblag}. Then it is clear (see for
instance lemma~1 in \cite{Soblag}) that $a_i(x)$ must be a polynomial in $x$
of degree at most $i$ for each $i=1,2,3,\ldots$. In \cite{Symjac} we showed
that the solution for $\l\{a_i(x)\r\}_{i=0}^{\infty}$ is not unique. In fact
it was shown that
\be\la{anul}a_0(x):=a_0(n,\a)=a_0(1,\a)b_0(n,\a)+c_0(n,\a),\;\ndots\ee
and that
\be\la{a}a_i(x)=a_0(1,\a)b_i(x)+c_i(x),\;i=1,2,3,\ldots,\ee
where $a_0(1,\a)$ is arbitrary and
\be\la{bnul}b_0(n,\a)=\frac{1}{2}\l[1-(-1)^n\r],\;\ndots,\ee
\be\la{cnul}c_0(n,\a)=4(2\a+3){n+2\a+2 \ch n-2},\;\ndots,\ee
\be\la{b}b_i(x)=\frac{2^{i-1}}{i!}(-x)^i,\;i=1,2,3,\ldots,\ee
\be\la{c}c_1(x)=0\;\mbox{ and }\;c_i(x)=(2\a+3)(1-x^2)
\frac{2^i}{i!}P_{i-2}^{(\a-i+3,\a-i+3)}(x),\;i=2,3,4,\ldots.\ee

In this paper we will give an alternative proof of this by using the
inversion formula (\ref{inv}).

By considering (\ref{sys1}) and (\ref{sys2}) for $n=0$ and $n=1$ we conclude
that $a_0(0,\a)=0$, $a_0(1,\a)$ is arbitrary and $a_1(x)=-a_0(1,\a)x$. For
$n=2,3,4,\ldots$ it turns out to be more convenient to use another system of
equations instead of (\ref{sys2}). By using (\ref{rel1}) we find for
$i=0,1,2,\ldots$
\bea & &iD^i\SP+xD^{i+1}\SP=D^i\l[xD\SP\r]\nn
&=&nD^i\SP+\frac{1}{2}(n+\a)D^iP_{n-2}^{(\a+1,\a+1)}(x),
\;n=2,3,4,\ldots.\n\eea
Combining (\ref{sys1}) and (\ref{sys2}) we now obtain
\be\la{sys3}\sum_{i=0}^{\infty}a_i(x)D^iP_{n-2}^{(\a+1,\a+1)}(x)=
\frac{8}{n+\a}{n+2\a \ch n-2}D^2\SP,\;n=2,3,4,\ldots.\ee
So we conclude that (\ref{sys2}) for $n=2,3,4,\ldots$ may be replaced by
(\ref{sys3}). Note that for $n=2$ this implies that $a_0(2,\a)=4(2\a+3)$.

Since $a_i(x)$ must be a polynomial in $x$ of degree at most $i$ for each
$i=1,2,3,\ldots$ we may write
$$a_i(x)=k_ix^i+\mbox{ lower order terms },\;i=1,2,3,\ldots.$$
By comparing the coefficients of highest degree in (\ref{sys1}) and
(\ref{sys3}) we find by using (\ref{def1})~:
$$\frac{a_0(n,\a)}{n!}+\sum_{i=1}^n\frac{k_i}{(n-i)!}=0,\;n=1,2,3,\ldots$$
and
$$\frac{a_0(n,\a)}{(n-2)!}+\sum_{i=1}^{n-2}\frac{k_i}{(n-i-2)!}=
4(2n+2\a-1){n+2\a \ch n-2}\frac{1}{(n-2)!},\;n=3,4,5,\ldots.$$
Since $k_i$ is independent of $n$ for $i=1,2,3,\ldots$ and
$a_0(2,\a)=4(2\a+3)$ we conclude that
\be\la{verschil}a_0(n+2,\a)-a_0(n,\a)=4(2n+2\a+3){n+2\a+2 \ch n},
\;\ndots,\ee
where $a_0(0,\a)=0$ and $a_0(1,\a)$ is arbitrary. Hence
$$a_0(2n,\a)-a_0(0,\a)=4\sum_{k=0}^{n-1}{2k+2\a+2 \ch 2k}(4k+2\a+3),
\;n=1,2,3,\ldots$$
and
$$a_0(2n+1,\a)-a_0(1,\a)=4\sum_{k=0}^{n-1}{2k+2\a+3 \ch 2k+1}(4k+2\a+5),
\;n=1,2,3,\ldots.$$
Note that we have
$$(2n+2\a+3){n+2\a+2 \ch n}=
(2\a+3)\l[{n+2\a+4 \ch n}-{n+2\a+2 \ch n-2}\r],\;\ndots.$$
Hence, by using the telescoping property of the sums we find that
$$\sum_{k=0}^{n-1}{2k+2\a+2 \ch 2k}(4k+2\a+3)=(2\a+3){2n+2\a+2 \ch 2n-2},
\;n=1,2,3,\ldots$$
and
$$\sum_{k=0}^{n-1}{2k+2\a+3 \ch 2k+1}(4k+2\a+5)=(2\a+3){2n+2\a+3 \ch 2n-1},
\;n=1,2,3,\ldots.$$
So we conclude that (\ref{anul}), (\ref{bnul}) and (\ref{cnul}) hold.

The systems of equations (\ref{sys1}) and (\ref{sys3}) lead to
\be\la{sys4}\sum_{i=1}^{\infty}a_i(x)D^i\SP=
\frac{4}{2\a+1}{n+2\a \ch n}D^2\SP-a_0(n,\a)\SP\ee
for $\ndots$ and
\be\la{sys5}\sum_{i=1}^{\infty}a_i(x)D^iP_{n-2}^{(\a+1,\a+1)}(x)=
\frac{8}{n+\a}{n+2\a \ch n-2}D^2\SP-a_0(n,\a)P_{n-2}^{(\a+1,\a+1)}(x)\ee
for $n=2,3,4,\ldots$.

First we remark that (\ref{sys4}) is true for $n=0$ and $n=1$ since
$a_0(0,\a)=0$ and $a_1(x)=-a_0(1,\a)x$. Then we will show that every
solution of (\ref{sys5}) also satisfies (\ref{sys4}). Suppose that
$\l\{a_i(x)\r\}_{i=1}^{\infty}$ is a solution of (\ref{sys5}). Now we use
(\ref{diff}), (\ref{rel2}), (\ref{verschil}) and the fact that
$\l\{a_i(x)\r\}_{i=1}^{\infty}$ are independent of $n$ to obtain for
$n=2,3,4,\ldots$ (see \cite{Report} for more details)
\bea & &2(n+\a+1)(2n+2\a+1)\sum_{i=1}^{\infty}a_i(x)D^i\SP\nn
&=&(n+2\a+1)(n+2\a+2)\sum_{i=1}^{\infty}a_i(x)D^iP_n^{(\a+1,\a+1)}(x)+{}\nn
& &{}\hspace{1cm}{}
-(n+\a)(n+\a+1)\sum_{i=1}^{\infty}a_i(x)D^iP_{n-2}^{(\a+1,\a+1)}(x)\nn
&=&(n+2\a+1)(n+2\a+2)\times{}\nn
& &{}\hspace{1cm}{}\times
\l[\frac{8}{n+\a+2}{n+2\a+2 \ch n}D^2P_{n+2}^{(\a,\a)}(x)
-a_0(n+2,\a)P_n^{(\a+1,\a+1)}(x)\r]+{}\nn
& &{}\hspace{1cm}{}-(n+\a)(n+\a+1)\times{}\nn
& &{}\hspace{3cm}{}\times\l[\frac{8}{n+\a}{n+2\a \ch n-2}D^2\SP
-a_0(n,\a)P_{n-2}^{(\a+1,\a+1)}(x)\r]\nn
&=&2(n+\a+1)(n+2\a+1)(n+2\a+2){n+2\a+2 \ch n}P_{n-2}^{(\a+2,\a+2)}(x)+{}\nn
& &{}\hspace{1cm}{}-8(n+\a+1){n+2\a \ch n-2}D^2\SP+{}\nn
& &{}\hspace{1cm}{}-2(n+\a+1)(2n+2\a+1)a_0(n,\a)\SP\nn
&=&2(n+\a+1)(2n+2\a+1)\l[\frac{4}{2\a+1}{n+2\a \ch n}D^2\SP
-a_0(n,\a)\SP\r].\n\eea
Since $\a>-1$ this proves that every solution of (\ref{sys5}) also satisfies
(\ref{sys4}).

Now we will solve (\ref{sys5}). Shifting $n$ by two we may write, since the
coefficients $\l\{a_i(x)\r\}_{i=1}^{\infty}$ are independent of $n$
\be\la{sys6}\sum_{i=1}^{\infty}a_i(x)D^iP_n^{(\a+1,\a+1)}(x)=F_n(x),
\;\ndots,\ee
where
$$F_n(x)=\frac{8}{n+\a+2}{n+2\a+2 \ch n}D^2P_{n+2}^{(\a,\a)}(x)
-a_0(n+2,\a)P_n^{(\a+1,\a+1)}(x).$$
Since $a_0(2,\a)=4(2\a+3)$ we easily find that $F_0(x)=0$. This implies that
the system of equations (\ref{sys6}) is of the form (\ref{sysalg}). So if we
apply the inversion formula (\ref{inv}) to the system of equations
(\ref{sys6}) we obtain by using (\ref{solalg})
$$a_i(x)=2^i\sum_{j=1}^i\frac{2\a+2j+3}{(2\a+j+3)_{i+1}}
P_{i-j}^{(-\a-i-2,-\a-i-2)}(x)F_j(x),\;i=1,2,3,\ldots.$$
Hence, by using (\ref{anul}) we conclude that the coefficients
$\l\{a_i(x)\r\}_{i=1}^{\infty}$ can be written in the form (\ref{a}).
Moreover, we find by using (\ref{bnul}), (\ref{sym}) and (\ref{spec})
\bea b_i(x)&=&2^{i-1}\sum_{j=1}^i\frac{2\a+2j+3}{(2\a+j+3)_{i+1}}
P_{i-j}^{(-\a-i-2,-\a-i-2)}(x)P_j^{(\a+1,\a+1)}(x)\l[(-1)^j-1\r]\nn
&=&2^{i-1}\sum_{j=0}^i\frac{2\a+2j+3}{(2\a+j+3)_{i+1}}
P_{i-j}^{(-\a-i-2,-\a-i-2)}(x)P_j^{(\a+1,\a+1)}(-x)+{}\nn
& &{}\hspace{1cm}{}-2^{i-1}\sum_{j=0}^i\frac{2\a+2j+3}{(2\a+j+3)_{i+1}}
P_{i-j}^{(-\a-i-2,-\a-i-2)}(x)P_j^{(\a+1,\a+1)}(x)\nn
&=&2^{i-1}\frac{(-x)^i}{i!},\;i=1,2,3,\ldots,\n\eea
which proves (\ref{b}). And by using (\ref{cnul}) and (\ref{diff}) we obtain
$$c_i(x)=2^i\sum_{j=1}^i\frac{2\a+2j+3}{(2\a+j+3)_{i+1}}
P_{i-j}^{(-\a-i-2,-\a-i-2)}(x)G_j(x),\;i=1,2,3,\ldots,$$
where
$$G_j(x)=\frac{4(2\a+3)}{j+\a+2}{j+2\a+4 \ch j}
\l[(\a+2)P_j^{(\a+2,\a+2)}(x)-(j+\a+2)P_j^{(\a+1,\a+1)}(x)\r].$$
It is clear that $G_1(x)=0$, which implies that $c_1(x)=0$. Note that since
$b_1(x)=-x$ this also implies that $a_1(x)=-a_0(1,\a)x$, which agrees with
what we have found before. Now we use (\ref{rel3}) to find
$$G_j(x)=(2\a+3)(1-x^2){j+2\a+4 \ch j}P_{j-2}^{(\a+3,\a+3)}(x),
\;j=2,3,4,\ldots.$$
Hence, for $i=2,3,4\ldots$ we have
\bea & &c_i(x)=(2\a+3)(1-x^2)2^i\times{}\nn
& &{}\hspace{2cm}{}\times\sum_{j=2}^i\frac{2\a+2j+3}{(2\a+j+3)_{i+1}}
{j+2\a+4 \ch j}P_{i-j}^{(-\a-i-2,-\a-i-2)}(x)P_{j-2}^{(\a+3,\a+3)}(x).\n\eea
Now it remains to show that
\bea\la{cc} & &\sum_{j=2}^i\frac{2\a+2j+3}{(2\a+j+3)_{i+1}}
{j+2\a+4 \ch j}P_{i-j}^{(-\a-i-2,-\a-i-2)}(x)P_{j-2}^{(\a+3,\a+3)}(x)\nn
& &{}\hspace{1cm}{}
=\frac{1}{i!}P_{i-2}^{(\a-i+3,\a-i+3)}(x),\;i=2,3,4,\ldots.\eea
In order to do this we write for $i=2,3,4,\ldots$
\bea & &\sum_{j=2}^i\frac{2\a+2j+3}{(2\a+j+3)_{i+1}}
{j+2\a+4 \ch j}P_{i-j}^{(-\a-i-2,-\a-i-2)}(x)P_{j-2}^{(\a+3,\a+3)}(x)\nn
& &{}\hspace{1cm}{}=\sum_{k=0}^{i-2}\frac{2\a+2k+7}{(2\a+k+5)_{i+1}}
\frac{(2\a+5)_{k+2}}{(k+2)!}P_{i-k-2}^{(-\a-i-2,-\a-i-2)}(x)
P_k^{(\a+3,\a+3)}(x).\n\eea
Now we apply definition (\ref{def1}) to $P_k^{(\a+3,\a+3)}(x)$ and
definition (\ref{def2}) to $P_{i-k-2}^{(-\a-i-2,-\a-i-2)}(x)$ and change
the order of summation to obtain for $i=2,3,4,\ldots$ (see \cite{Report})
\bea & &\sum_{k=0}^{i-2}\frac{2\a+2k+7}{(2\a+k+5)_{i+1}}
\frac{(2\a+5)_{k+2}}{(k+2)!}P_{i-k-2}^{(-\a-i-2,-\a-i-2)}(x)
P_k^{(\a+3,\a+3)}(x)\nn
&=&\sum_{k=0}^{i-2}\sum_{m=0}^{i-k-2}\sum_{n=0}^k(-1)^{i-k-2}\times{}\nn
& &{}\hspace{1cm}{}\times\frac{(2\a+2k+7)(2\a+5)_{k+n+2}(\a+n+4)_{i-m-n-2}}
{(2\a+k+5)_{i-m+1}(k+2)!\,m!\,(i-k-m-2)!\,n!\,(k-n)!}\l(\frac{x-1}{2}\r)^{m+n}\nn
&=&\sum_{j=0}^{i-2}\sum_{n=0}^j(-1)^{i-n}\frac{(2\a+5)_n(\a+n+4)_{i-j-2}}
{(2\a+2n+7)_{i-j-1}(n+2)!\,(j-n)!\,(i-j-2)!\,n!}\l(\frac{x-1}{2}\r)^j\times{}\nn
& &{}\hspace{1cm}{}\times\sum_{k=0}^{i-j-2}
\frac{(-i+j+2)_k(2\a+n+5)_k(2\a+2n+7)_k}{(2\a+2n+i-j+6)_k(n+3)_kk!}
(2\a+2n+2k+7).\n\eea
In \cite{Report} we proved that for $i-j-2\in\{0,1,2,\ldots\}$ we have
\bea & &\sum_{k=0}^{i-j-2}
\frac{(-i+j+2)_k(2\a+n+5)_k(2\a+2n+7)_k}{(2\a+2n+i-j+6)_k(n+3)_kk!}
(2\a+2n+2k+7)\nn
& &{}\hspace{1cm}{}=\frac{(2\a+2n+7)_{i-j-1}(-\a-1)_{i-j-2}}
{(n+3)_{i-j-2}(\a+n+4)_{i-j-2}}.\n\eea
See \cite{Report} for much more details. By using this, the well-known
Vandermonde summation formula and definition (\ref{def2}) we finally obtain
for $i=2,3,4,\ldots$
\bea & &\sum_{k=0}^{i-2}\frac{2\a+2k+7}{(2\a+k+5)_{i+1}}
\frac{(2\a+5)_{k+2}}{(k+2)!}P_{i-k-2}^{(-\a-i-2,-\a-i-2)}(x)
P_k^{(\a+3,\a+3)}(x)\nn
&=&\sum_{j=0}^{i-2}\sum_{n=0}^j(-1)^{i-n}\frac{(2\a+5)_n(-\a-1)_{i-j-2}}
{(n+i-j)!\,(j-n)!\,(i-j-2)!\,n!}\l(\frac{x-1}{2}\r)^j\nn
&=&(-1)^i\sum_{j=0}^{i-2}\frac{(-\a-1)_{i-j-2}}{(i-j)!\,(i-j-2)!\,j!}
\l(\frac{x-1}{2}\r)^j\hyp{2}{1}{-j,2\a+5}{i-j+1}{1}\nn
&=&\frac{(-1)^i}{i!}\sum_{j=0}^{i-2}
\frac{(-\a-1)_{i-j-2}(i-j-2\a-4)_j}{(i-j-2)!\,j!}\l(\frac{x-1}{2}\r)^j\nn
&=&\frac{(-1)^{i-2}}{i!}\sum_{j=0}^{i-2}\frac{(i-j-2\a-4)_j}{j!}
\frac{(-\a-1)_{i-j-2}}{(i-j-2)!}\l(\frac{x-1}{2}\r)^j
=\frac{1}{i!}P_{i-2}^{(\a-i+3,\a-i+3)}(x),\n\eea
which proves (\ref{cc}).

\section{Some remarks}

By using definition (\ref{def3}) we may write
$$P_{i-2}^{(\a-i+3,\a-i+3)}(x)=2^{2-i}\sum_{k=0}^{i-2}{\a+1 \ch i-2-k}
{\a+1 \ch k}(x-1)^k(x+1)^{i-2-k},\;i=2,3,4,\ldots.$$
By using (\ref{c}) this implies that for nonnegative integer values of $\a$
we have $c_i(x)=0$ for all $i>2\a+4$ and
$$c_{2\a+4}(x)=(2\a+3)(1-x^2)\frac{2^{2\a+4}}{(2\a+4)!}
P_{2\a+2}^{(-\a-1,-\a-1)}(x)=-\frac{4(2\a+3)}{(2\a+4)!}(x^2-1)^{\a+2}\ne 0.$$
Hence, for nonnegative integer values of $\a$ the polynomials $\set{\SGP}$
satisfy a differential equation of the form
$$M\sum_{i=0}^{2\a+4}c_i(x)y^{(i)}(x)+
(1-x^2)y''(x)-2(\a+1)xy'(x)+n(n+2\a+1)y(x)=0,$$
which is of finite order $2\a+4$ if $M>0$. Moreover, in \cite{Symjac} it was
shown that for $M>0$ the differential equation given by (\ref{DVSGP}) has
finite order if and only if we choose $a_0(1,\a)=0$ and if $\a$ is a
nonnegative integer.

Finally, we remark that if we apply the inversion formula (\ref{inv}) to
the system of equations (\ref{sys4}) instead of (\ref{sys5}) we find for
$i=1,2,3,\ldots$ that
\be\la{b*}b_i(x)=2^{i-1}\sum_{j=1}^i\frac{2\a+2j+1}{(2\a+j+1)_{i+1}}
P_{i-j}^{(-\a-i-1,-\a-i-1)}(x)P_j^{(\a,\a)}(x)\l[(-1)^j-1\r]\ee
and
\bea\la{c*} & &c_i(x)=2^{i+2}\sum_{j=1}^i\frac{2\a+2j+1}{(2\a+j+1)_{i+1}}
P_{i-j}^{(-\a-i-1,-\a-i-1)}(x)\times{}\nn
& &{}\hspace{1cm}{}\times\l[\frac{1}{2\a+1}{j+2\a \ch j}D^2P_j^{(\a,\a)}(x)
-(2\a+3){j+2\a+2 \ch j-2}P_j^{(\a,\a)}(x)\r].\eea
From (\ref{b*}) we easily obtain (\ref{b}) in the same way as before by
using (\ref{sym}), (\ref{inv}) and (\ref{spec}). Further we easily
find from (\ref{c*}) that $c_1(x)=0$, but we were not able to derive
(\ref{c}) for $i=2,3,4,\ldots$ from (\ref{c*}).



\begin{flushleft}
{ J.~Koekoek } \\
{ Menelaoslaan 4, 5631 LN Eindhoven, The Netherlands }  \\
{ }
\end{flushleft}

\begin{flushleft}
{ R.~Koekoek } \\
{ Delft University of Technology, Faculty of Technical Mathematics and
Informatics, }\\
{ P.O. Box 5031, 2600 GA Delft, The Netherlands }\\
{ e-mail : koekoek@twi.tudelft.nl }
\end{flushleft}

\end{document}